\documentclass[12pt,notitlepage,twoside]{article}
\usepackage{latexsym} 
\usepackage{amsmath} 
\usepackage{amsfonts} 
\usepackage{amssymb} 
\renewcommand{\a }{\alpha } 
\renewcommand{\b }{\beta } 
\renewcommand{\d}{\delta } 
\newcommand{\media}{\mkern12mu\hbox{\vrule height4pt depth-3.2pt width5pt} \mkern-16mu\int}
\newcommand{\laplconf}{L_{g}}
\newcommand{\qi}{q_i}
\newcommand{\yi}{y_i}
\newcommand{\vi}{v_i}
\newcommand{\wi}{w_i}
\newcommand{\ui}{u_i}
\newcommand{\D }{\Delta }
 
\newcommand{\e }{\varepsilon }

\newcommand{\graffe}[1]{\parent \{ \}{#1}} 
 \newcommand{\foral }{\forall\, } 
\newcommand{\G }{\Gamma } 
\renewcommand{\l }{\lambda }

\newcommand{\n }{\nabla } 
\newcommand{\var }{\varphi } 
 
\newcommand{\s }{\sigma }

\newcommand{\ov}{\overline} 
\newcommand{\wtilde }{\widetilde} 
 
\newcommand{\be}{\begin{equation}} 
\newcommand{\ee}{\end{equation}} 
\newenvironment{pf}{\noindent{\bf Proof.}\enspace}{
\hfill$\Box$\medskip} 
\newenvironment{pfn}[1]{\noindent{\bf Proof of {#1}\enspace}}{
\hfill$\Box$\medskip} 
\newcommand{\R}{\mathbb{R}} 
 
 \newcommand{\parent}[3]{\left #1 {#3} \right #2} 
 \newcommand{\tonde}[1]{\parent (){#1}} 
 \newcommand{\quadre}[1]{\parent []{#1}} 
\newcommand{\barre}[1]{\parent \Vert \Vert {#1}} 

\newcommand{\N}{\mathbb{N}}

\newtheorem{thm}{Theorem}[section] 
\newtheorem{pro}[thm]{Proposition}
\newtheorem{lem}[thm]{Lemma}

\newtheorem{cor}[thm]{Corollary}
\newtheorem{df}[thm]{Definition} 
\numberwithin{equation}{section}

\author{{\sc  Veronica Felli and Mohameden Ould Ahmedou}}

\title { \Large \textbf{On a geometric equation with critical\\
 nonlinearity on the boundary\\
} }
 
\begin{document}

\date{}

\maketitle

 
{\footnotesize 
\begin{abstract} 
 
\noindent 
A theorem of Escobar asserts that, on a positive three dimensional
smooth compact Riemannian manifold with boundary which is not
conformally equivalent to the standard three dimensional ball, a
necessary and sufficient condition for a $C^2$ function $H$ to be the
mean curvature of some conformal scalar flat metric is that $H$ is positive
somewhere. We show that, when the boundary is umbilic and the function $H$
is positive everywhere, all such metrics stay in a compact set with respect to the $C^2$ norm and the total degree of all solutions is equal to $-1$.


\medskip\noindent\footnotesize {{\bf MSC classification:}\quad 35J60, 53C21, 58G30.}

\end{abstract}

}

\section{Introduction}
\mbox{}

In \cite{E1}, Jos\'e F. Escobar raised the following question:
given a compact Riemannian manifold with boundary, when is it conformally equivalent to one that has zero scalar curvature and whose boundary has a constant mean curvature? This problem can be seen as a ``generalization'' to higher dimensions of the well known Riemannian Mapping Theorem. The later states that an open, simply connected proper subset of the plane is conformally diffeomorphic to the disk. In higher dimensions few regions are conformally diffeomorphic to the ball. However one can still ask whether a domain is conformal to a manifold that resembles  the ball in two ways: namely, it has zero scalar curvature and its boundary has constant mean curvature. The above problem is equivalent to finding a smooth positive solution $u$ to the following nonlinear boundary value problem on a Riemannian manifold with boundary $(M^n,g)$, $n \geq 3$: 
\begin{equation}\tag{$P$}
\begin{cases}
- \Delta_g u + \displaystyle{\frac{(n-2)}{4(n-1)}}R_g u = 0, \quad u > 0,  & 
\mbox{in}\ \mathring{M}, \\[10pt]
\displaystyle{\frac{\partial u}{\partial\nu}}  + \displaystyle{\frac{n - 2}{2}} h_g u = c u^{\frac{n}{n-2}}, & 
\mbox{on}\ \partial M,  
\end{cases}
\end{equation}
 where $R_g$ is the scalar curvature of $M$, $h_g$ is the mean curvature of $ \partial M$, $\nu$ is the outer normal vector with respect to $g$, and $c$ is a constant whose sign is uniquely determined by the conformal structure.

For almost all manifolds, Escobar \cite{E1,E2} established that $(P)$ has a solution. More recently in \cite{Ahm} this problem has been studied using the tools of the critical points at infinity of A. Bahri \cite{abb}, see also Bahri-Coron \cite{bc} and Bahri-Brezis \cite{bb}. Going beyond the existence results of the above paper, we proved recently in \cite{VA} that, when $(M,g)$ is locally conformally flat with umbilic boundary but not conformal to the standard ball, all solutions of $(P)$ stay in a compact set with respect to the $C^2$ norm and the total degree of all solutions is equal to $-1$.

The heart of the proof of the above result is some fine analysis of possible blow-up behaviour of solutions to $(P)$. More specifically we obtained energy independent estimates of solutions to 
\[
  \begin{cases} 
   L_gu=0,\quad u>0,         &\textrm{in}\  \mathring{M},\\
    B_gu=(n-2)u^q,            &\textrm{on}\ \partial M,
   \end{cases}
\]
where 
\[
  1<1+\e_0\leq q\leq \frac {n}{n-2},\quad L_g=\D_g-\frac{n-2}{4(n-1)}R_g, \quad B_g=\frac{\partial}{\partial \nu_g}+\frac{n-2}{2}h_g.
\]
Instead of looking for conformal metrics with zero scalar curvature
and constant mean curvature as in $(P)$, one may also look for scalar
flat conformal metrics with boundary mean curvature being a given
function $H$; this problem is equivalent to finding a smooth positive
solution $u$ to
\begin{equation}\tag{$P_H$}
\begin{cases}
L_gu=0,\quad u>0,     &\textrm{in}\ \mathring M,\\
B_gu=Hu^{\frac n{n-2}},            &\textrm{on}\ \partial M.
\end{cases}
\end{equation}
Such a problem was studied in \cite{E2} by Escobar, who proved that on
a positive three dimensional smooth compact Riemannian manifold which
is not conformally equivalent to the standard three dimensional ball,
a necessary and sufficient condition for a $C^2$ function $H$ to be
the mean curvature of some conformal flat metric is that $H$ is
positive somewhere. We recall that a manifold is called of positive
type if the quadratic part of the Euler functional associated to $(P)$
is positive definite. \\
In our work we assume that the boundary is
umbilic, that is the traceless part of the second fundamental form
vanishes on the boundary. Moreover we assume that the function $H$ is
positive.

Our first theorem gives a priori estimates of solutions of $(P_{H,q})$ in $H^1(M)$ norm.
\begin{thm}\label{t:sth1}
Let $(M,g)$ be a three dimensional smooth compact Riemannian manifold
with umbilic boundary. Then for all $\e_0>0$
\[
 \|u\|_{H^1(M)}\leq C\quad\foral u\in\bigcup_{1+\e_o\leq q\leq3}{\cal M}_{H, q},
\]
where $C$ depends only on $M$, $g$, $\e_0$, $\|H\|_{C^2(\partial M)}$, and the positive lower bound of~$H$.
\end{thm}

Our next theorem states that for any positive $C^2$ function $H$, all such metrics stay bounded with respect to the $C^2$ norm and the total Leray-Schauder degree of all the solutions of $(P_H)$ is $-1$. In fact we establish a slightly stronger compactness result. Consider for $1<q\leq 3$ the problem
\begin{equation}\tag{$P_{H,q}$}
\begin{cases}
L_gu=0,\quad u>0,     &\textrm{in $\mathring M$,}\\
B_gu=Hu^q,            &\textrm{on $\partial M$.}
\end{cases}
\end{equation}
We use ${\cal M}_{H, q}$ to denote the set of solutions of $P_{H,q}$
in $C^2(M)$. We have the following theorem. 
\begin{thm}\label{t:pos}
Let $(M,g)$ be a positive three dimensional smooth compact Riemannian
manifold with umbilic boundary which is not conformally equivalent to the standard three dimensional ball. Then, for any $1<q\leq 3$ and positive function $H\in C^2(\partial M)$, there exists some constant $C$ depending only on $M,g,\|H\|_{C^2}$, the positive lower bound of $H$ and $q-1$ such that
\[
  \frac1C\leq u\leq C \quad\textrm{and}\quad \|u\|_{C^2(M)}\leq C
\]
for all solutions $u$ of $(P_{H,q})$. Moreover the total degree of all solutions of $(P_{H,q})$ is $-1$. Consequently, equation $(P_{H,3})$ has at least one solution.
\end{thm}
We remark that the hypothesis that $(M,g)$ is not conformally
equivalent to the standard three dimensional ball is necessary since
$(P_H)$ may have no solution in this case due to the
Kazdan-Warner's conditions for solvability. On the ball sufficient
conditions on $H$ in dimensions $3$ and $4$ are given in \cite{DMOA2}
and \cite{EG}, and perturbative results were obtained in \cite{CXY}.

Finally, let us point out that recently S. Brendle \cite{brendle1,brendle2} obtained on surfaces some results related to ours. He used curvature flows methods, in the spirit of M. Struwe \cite{S} and X. X. Chen \cite{C}. The curvature flow method was introduced by R. Hamilton \cite{H}, and used by B. Chow \cite{chow}, R. Ye \cite{ye}, and Bartz-Struwe-Ye \cite{BSY}.

The remainder of the paper is organized as follows. In section 2 we provide the main local blow-up analysis giving first sharp pointwise estimates to a sequence of solutions near isolated simple blow-up points, then we prove that an isolated blow-up is in fact an isolated simple blow up, ruling out the possibility of bubbles on top of bubbles. In section~3 we rule out the possibility of bubble accumulations and establish Theorem \ref{t:sth1}. In section 4 we study compactness of solutions of $(P_H)$ and establish Theorem \ref{t:pos}. In the Appendix, we provide some standard descriptions of singular behaviour of positive solutions to some linear boundary value elliptic equations in punctured half balls and collect some useful results. 

\section{Local blow-up analysis}
\mbox{}

We may assume without loss of generality that $h_g\equiv0$. Indeed, let $\var_1$ be a positive eigenfunction associated to the first eigenvalue $\l_1$ of the problem 
\[
\begin{cases}
L_g\var=\l_1\varphi,&\text{in}\ \mathring M,\\
B_g\var=0,&\text{on}\ \partial M.
\end{cases}
\]
Setting $\tilde g=\var_1^4g$ and $\tilde u=\var_1^{-1}u$, where $u$ is
a solution of $(P_{H,3})$, one can easily check that $R_{\tilde g}>0$,
$h_{\tilde g}\equiv0$, and $\tilde u$ satisfies
\[
\begin{cases}
L_{\tilde g}\tilde u =0,&\text{in}\ \mathring M,\\
\displaystyle{\frac{\partial\tilde u}{\partial\nu}}=H\tilde u^3,&\text{on}\ \partial M.
\end{cases}
\]
For sake of simplicity, we work with $\tilde g$, denoting it by
$g$. Since $\partial M$ is umbilic with respect to $g$, and $h_{\tilde
  g}=0$, it follows that the second fundamental form vanishes at each
point of the boundary, that is the boundary is a totally geodesic
submanifold. Hence we can take conformal normal coordinates around any point
of the boundary \cite{E3}

Let us first recall the definitions of isolated and isolated simple blow up which were first introduced by R. Schoen \cite{sch} and used extensively by Y. Y. Li \cite{liuno,lidue}.

\begin{df}\label{d:isolated}
Let $(M,g)$ be a smooth compact $n$-dimensional Riemannian manifold with boundary and let $\bar r>0$, $\bar c>0$, $\bar x\in\partial M$, $H\in C^0(\overline{B_{\bar r}(\bar x)})$ be some positive function where  $B_{\bar r}(\bar x)$ denotes the geodesic ball in $(M,g)$ of radius $\bar r$ centered at $\bar x$. Suppose that, for some sequences $q_i=3-\tau_i$, $\tau_i\to 0$, $H_i\to H$ in $C^2(\overline{B_{\bar r}(\bar x)})$, $\{u_i\}_{i\in \N}$ solves
\begin{equation}\label{e:seq}
\begin{cases}
\laplconf u_i=0,\quad\ui>0, &\textrm{in $B_{\bar r}(\bar x)$,}\\
 \displaystyle{\frac{\partial u_i}{\partial \nu}}=H_i\ui^{\qi}, &\textrm{on $\partial M\cap B_{\bar r}(\bar x)$.}
\end{cases}
\end{equation} 
We say that ${\bar x}$ is an isolated blow-up point of $\{\ui\}_i$ if there exists a sequence of local maximum points $x_i$ of $\ui$ such that $x_i\to \bar x$ and, for some $C_1>0$,
\[
\lim_{i\to\infty}\ui(x_i)=+\infty\quad\textrm{and}\quad\ui(x)\leq  C_1 d(x,x_i)^{-\frac{1}{q_i-1}},\quad\foral x\in B_{\bar r}(x_i),\ \foral i.
\]
\end{df}
To describe the behaviour of blowing-up solutions near an isolated blow-up point, we define spherical averages of $\ui$ centered at $x_i$ as follows
\[
\bar \ui(r)=\media_{M\cap\partial B_r(\bar x)}\ui=\frac{1}{{\rm Vol}_g(M\cap\partial B_r(\bar x))}\int_{M\cap\partial B_r(\bar x)}\ui.
\]
Now we define the notion of isolated simple blow-up point.
\begin{df}\label{d:isolatedsimple}
Let $x_i\to \bar x$ be an isolated blow-up point of $\{\ui\}_i$ as in Definition \ref{d:isolated}. We say that $x_i\to \bar x$ is an isolated simple blow-up point of $\{\ui\}_i$ if, for some positive constants $\tilde r\in(0,\bar r)$ and $C_2>1$, the function $\bar w_i(r):=r^{\frac1{q_i-1}}\bar\ui(r)$ satisfies, for large $i$,
\[
\bar w_i'(r)<0\quad\textrm{for $r$ satisfying $C_2\ui^{1-q_i}(x_i)\leq r\leq\tilde r$.}
\]
\end{df}
For any $\bar x\in\partial M$, by choosing geodesic normal coordinate system centered at $\bar x$, we can assume without loss of generality that 
$$
\displaylines{\bar x=0, \quad g_{ij}(0)=\d_{ij},\quad B_1^+(0):=\{x=(x^1,x^2,x^3):\ |x|<1\ \mbox{and}\ x^3>0\}\subset M,\cr
\{(x',0)=(x^1,x^2,0):\ |x'|<1\}\subset \partial M,\quad\G_{ij}^k(0)=0,\cr}
$$
where $\G_{ij}^k$ is the Christofell symbol.
For later use, we denote 
\begin{eqnarray*}
&\R^3_+=\{(x',x^3)\in \R^2\times \R:\ x^3>0\},\quad B^+_r(\bar x)=\{x=(x',x^3)\in\R^3_+:\ |x-\bar x|<r\},\\
&B^+_r=B_r^+(0),\quad \Gamma_1(B_r^+(\bar x))=\partial B_r^+(\bar x)\cap\partial\R^3_+,\quad \Gamma_2(B_r(\bar x))=\partial B_r(\bar x)\cap\R^3_+.
\end{eqnarray*}
Let $H_i\to H$ in $C^2(\G_1(B_3^+))$ be a sequence of positive functions, $\qi$ be a sequence of numbers satisfying $2\leq q_i\leq 3$ and $\qi\to3$, and $\{\vi\}_i\subset C^2(\overline{B_3^+})$ be a sequence of solutions to
\begin{equation}\tag{$P_i$}
\begin{cases}
\displaystyle{-\Delta_g\vi+\frac 18 R_g v_i=0,\quad \vi>0}, &\textrm{in $B_3^+$,}\\
\displaystyle{\frac{\partial\vi}{\partial\nu}}=H_i\vi^{\qi},&\textrm{on $\G_1(B_3^+)$.}
\end{cases}
\end{equation}
In this section, we start giving some properties of isolated and
isolated simple blow-up.  Hence forward we use~$c$ to denote positive
constants which may vary from formula to formula and which may depend
only on $M$, $g$, and $\bar r$. Such blow-up analysis was also carried
out in \cite{EG}, where $(M,g)$ was the standard ball endowed with
euclidean metric, see also our previous work \cite{VA}.

The following lemma gives a Harnack Inequality, which proof is
contained in~\cite{VA, EG}, Lemma 2.3, up to some minor modifications.
\begin{lem}\label{l:harnack}
Let $\vi$ satisfy $(P_i)$ and $\yi \to\bar y\in \G_1(B_3^+)$ be an isolated blow-up of $\{\vi\}_i$. Then for any $0<r<\bar r$, we have 
\[
\max_{\ov{B^+_{2r}(y_i)}\setminus B^+_{r/2}(y_i)}\vi\leq C_3\min_{\ov{B^+_{2r}(y_i)}\setminus B^+_{r/2}(y_i)}\vi,
\]
where $C_3$ is some positive constant independent of $i$ and $r$.
\end{lem}
\begin{lem}\label{l:bubble}
Let $\vi$ satisfy $(P_i)$, $\yi\to \bar y\in\G_1(B_1^+)$ be an isolated blow-up point. Then for any $R_i\to+\infty$, $\e_i\to 0^+$ we have that, after passing to a subsequence,
\begin{align}\label{e:bubble}
&\barre{\vi^{-1}(\yi)\vi\big(\exp_{\yi}(\vi^{1-\qi}(\yi)x)\big)-\displaystyle{\quadre{\frac1{(1+h_ix^3)^2+h_i^2|x'|^2}}}^{1/2}}_{C^1(B^+_{2R_i})}\notag\\
&\quad+\barre{\vi^{-1}(\yi)\vi\big(\exp_{\yi}(\vi^{1-\qi}(\yi)x)\big)-\displaystyle{\quadre{\frac1{(1+h_ix^3)^2+h_i^2|x'|^2}}}^{1/2}}_{H^1(B^+_{2R_i})}\leq\e_i
\end{align}
and
\begin{equation}\label{e:raggi}
\frac {R_i}{\log\vi(\yi)}\ \mathop{\longrightarrow}\limits_{i\to+\infty}\ 0,
\end{equation}
where $x=(x',x^3)\in B_1^+$ and $h_i~=~H_i(\yi)$.
\end{lem}
\begin{pf} Let $g_i=(g_i)_{\a\b}(x)\,dx^\a dx^\b=g_{\a\b}(\vi^{1-\qi}(\yi)x)\,dx^\a dx^\b$ denote the scaled metric. Set 
\[
\xi_i(x)=\vi^{-1}(\yi)\vi\tonde{\yi+\vi^{1-\qi}(\yi)x},\quad\textrm{for $x\in B^{-T_i}_{\vi^{\qi-1}(\yi)}$,}
\]
defined on the set
\[
B^{-T_i}_{\vi^{\qi-1}(y_i)}:=\left\{z\in\R^3:\ |z|<\vi^{\qi-1}(y_i)\quad\textrm{and}\quad z^3>-T_i\right\}
\]
where $T_i=y_i^3\vi^{\qi-1}(y_i)$. Then $\xi_i(x)$ satisfies
\begin{equation}\label{e:xi}
\begin{cases}
\displaystyle{-\Delta_{g_i}\xi_i+\frac 18\, \vi^{2(1-\qi)}(\yi)R_{g_i}(\yi+\vi^{1-\qi}(\yi)x)\xi_i=0},\quad\xi_i>0, \quad\textrm{in $B^{-T_i}_{\vi^{\qi-1}(y_i)}$,}\\
\displaystyle{\frac{\partial\xi_i} {\partial\nu_{g_i}}}=H_i(\yi+\vi^{1-\qi}(\yi)x)\xi_i^{\qi}\quad\textrm{on $\partial B^{-T_i}_{\vi^{\qi-1}(y_i)}\cap\{z\in\R^3:\ z^3=-T_i\}$,}\\
\xi_i(0)=1,\\
\textrm{$0$ is a local maximum point of $\xi_i$,}\\
0<\xi_i(x)\leq \tilde c|x|^{-\frac1{\qi-1}},
\end{cases}
\end{equation}
for some positive constant $\tilde c$. Now we prove that $\xi_i$ is
locally bounded. Using Hopf Point Boundary Lemma and Lemma~\ref{l:harnack}, we derive that for $0<r<1$
\[
1=\xi_i(0)\geq\min_{\ov{\G_1(B^+_r)}}\xi_i\geq\min_{\ov{\G_2(B^+_r)}}\xi_i\geq c\max_{\ov{\G_2(B^+_r)}}\xi_i
\]
which implies that, for some $c$ independent of $r$,
\[
\max_{\ov{\G_2(B^+_r)}}\xi_i\leq c.
\]
Therefore, we derive easily that $\xi_i$ is locally bounded. Applying standard elliptic estimates to $\{\xi_i\}$, we have, after passing to a subsequence, that $\xi_i\to\xi$ in $C^2_{\rm loc}(\R^3_+)$ and $H^1_{\rm loc}(\R^3_+)$ for some $\xi$ satisfying
\[
\begin{cases}
\Delta \xi=0,\quad\xi>0,&\textrm{in $\R^3_{-T}$,}\\
\displaystyle{\frac{\partial\xi}{\partial \nu}}=\quadre{\lim_iH_i(\yi)}\xi^3,&\textrm{on $\partial \R^3_{-T}$,}
\end{cases}
\]
where $\R^3_{-T}:=\{x=(x',x^3)\in\R^3:\ x^3>-T\}$ and $T=\lim_iT_i$. By the Liouville Theorem and the last estimate of \eqref{e:xi} we have that $T<+\infty$. By Li-Zhu \cite{LZ1} Liouville type Theorem (see Theorem \ref{t:liouville} of the Appendix), we easily deduce that $T=0$ and
\[
\xi(x',x^3)=
\quadre{\frac1{(1+\lim_iH_i(\yi)x^3)^2+(\lim_iH_i(\yi))^2|x'|^2}}^{1/2}.
\]
Lemma \ref{l:bubble} is proved.
\end{pf}

Before stating our next result, we point out that it follows from Lemma \ref{l:G} of the Appendix that, for $\d_0>0$ small enough, there exists a unique function $G(\cdot,\bar y)\in C^2(\ov{B^+_{\d_o}(\bar y)}\setminus\{\bar y\})$ satisfying
\[
\begin{cases}
\displaystyle{-\Delta_g G(\cdot,\bar y)+\frac18\,R_g G(\cdot,\bar y) =0},&\textrm{in $B^+_{\d_o}(\bar y)$,}\\
\displaystyle{\frac{\partial}{\partial\nu}}G(\cdot,\bar y)=0,&\textrm{on $\G_1(B^+_{\d_o}(\bar y))\setminus\{\bar y\}$,}\\
\lim_{y\to\bar y}d(y,\bar y)G(y,\bar y)=1.
\end{cases}
\]
Now we state our main estimate on isolated simple blow-up points.
\begin{pro}\label{p:upperbound}
Let $\vi$ satisfy $(P_i)$ and $\yi\to\bar y\in\G_1(B_1^+)$ be an isolated simple blow-up point, with \eqref{e:bubble} and \eqref{e:raggi} for all $i$. Then for some positive constant $C$ depending only on $C_1$, $\tilde r$, $\|H_i\|_{C^2(\G_1(B^+_3))}$, and $\inf_{y\in\G_1(B^+_1)}H_i(y)$  we have
\begin{equation}\label{e:upperbound}
\vi(y)\leq C\vi^{-1}(\yi)d(y,\yi)^{-1},\quad\textrm{for}\quad d(y,\yi)\leq\frac{\tilde r}2
\end{equation}
where $C_1$ and $\tilde r$ are given in Definitions \ref{d:isolated} and \ref{d:isolatedsimple}. Furthermore, after passing to some subsequence, for some positive constant $b$,
\[
\vi(\yi)\vi\ \mathop{\longrightarrow}\limits_{i\to+\infty}\ b\, G(\cdot,\bar y)+E\quad\textrm{in}\  C^2_{\rm loc}(\ov{B^+_{\tilde \rho}(\bar y)}\setminus\{\bar y\})
\]
where $\tilde\rho=\min(\delta_0,\tilde r/2)$ and $E\in C^2(B^+_{\tilde\rho}(\bar y))$ satisfies
\[
\begin{cases}
\displaystyle{-\Delta_g E+\frac 18\, R_g E=0},&\textrm{in $B^+_{\tilde\rho}$,}\\
\displaystyle{\frac{\partial E}{\partial\nu}}=0,&\textrm{on $\G_1(B^+_{\tilde\rho})$.}
\end{cases}
\]
\end{pro}
Proposition \ref{p:upperbound} will be established through a series of lemmas.
\begin{lem}\label{l:bounds}
Let $\vi$ satisfy $(P_i)$ and $\yi\to\bar y\in\G_1(B_1^+)$ be an isolated simple blow-up. Assume $R_i\to+\infty$ and $0<\e_i<e^{-R_i}$ are sequences for which \eqref{e:bubble} and \eqref{e:raggi} hold. Then for any given $0<\d<1/100$, there exists $\rho_1\in (0,\tilde r)$ which is independent of $i$ (but depending on $\delta$), such that 
\begin{align}
\vi(\yi)&\leq C_4 \vi^{-\l_i}(\yi)d(y,\yi)^{-1+\d},  &\foral r_i\leq d(y,\yi)\leq\rho_1,\label{e:funcbound}\\
 \n_g \vi(\yi)&\leq C_4 \vi^{-\l_i}(\yi)d(y,\yi)^{-2+\d},  &\foral r_i\leq d(y,\yi)\leq\rho_1,\label{e:gradbound}\\ 
 \n^2_g \vi(\yi)&\leq C_4 \vi^{-\l_i}(\yi)d(y,\yi)^{-3+\d},  &\foral r_i\leq d(y,\yi)\leq\rho_1,\label{e:grad2bound}
\end{align}
where $r_i=R_i\vi^{1-\qi}(\yi)$, $\l_i=(1-\d)(\qi-1)-1$, and $C_4$ is some positive constant independent of $i$.
\end{lem}
\begin{pf}
We assume, for simplicity, that $g$ is the flat metric. The general case can be derived essentially in the same way. Let $r_i=R_i\vi^{1-\qi}(\yi)$, it follows from Lemma \ref{l:bubble} that
\begin{equation}\label{e:efb}
\vi(y)\leq c\vi(\yi)R_i^{-1},\quad\textrm{for}\ d(y,\yi)=r_i.
\end{equation}
We then derive from Lemma \ref{l:harnack}, \eqref{e:efb}, and the definition of isolated simple blow-up that, for $r_i\leq d(y,\yi)\leq\tilde r$, we have
\begin{equation}\label{e:oe}
\vi^{\qi-1} (y)\leq cR_i^{-1+o(1)}d(y,\yi)^{-1}.
\end{equation}
Set $T_i=y_i^3\vi^{\qi-1}(\yi)$. From the proof of Lemma \ref{l:bubble} we know that $\lim_i T_i=0$. It is not restrictive to take $\yi=(0,0,y_i^3)$. Thus we have $d(0,y_i^3)=o(r_i)$. So 
\[
B_1^+(0)\setminus B_{2r_i}^+(0)\subset \graffe{\frac32 r_i\leq d(y,\yi)\leq \frac 32}.
\]
Let us apply the Maximum Principle stated in Theorem \ref{t:maxprin} in the Appendix; to this aim set 
\[
 \var_i(y)=M_i\tonde{|y|^{-\d}-\e|y|^{\d-1}y^3}+A\vi^{-\l_i}(\yi)\tonde{|y|^{-1+\d}-\e|y|^{-2+\d}y^3}
\]
with $M_i$ and $A$ to be chosen later, and let $\Phi_i$ be the boundary operator defined by
\[
\Phi_i(v)=\displaystyle{\frac{\partial v}{\partial\nu}}-H_i \vi^{\qi-1}(\yi)v.
\]
A direct computation yields
\[\Delta \var_i(y)=M_i|y|^{-\d}[-\d(1-\d)+O(\e)]+|y|^{-(3-\d)}A\vi^{-\l_i}(\yi)[-\d(1-\d)+O(\e)].
\]
So one can choose $\e=O(\d)$ such that $\D\var_i\leq0$.

Another straightforward computation taking into account \eqref{e:oe} shows that for $\d>0$ there exists $\rho_1(\d)>0$ such that
\[
\Phi_i\var_i>0\quad\textrm{on}\ \G_1(B^+_{\rho_1}).
\]
Taking 
\begin{align*}
\Omega&=D_i=B^+_{\rho_1}\setminus B^+_{2r_i}(0)&&{}\\
\Sigma&=\G_1(D_i),&&\G=\G_2(D_i),\\
V&\equiv0,&&h=H_i\vi^{\qi-1},\\
\psi&=\vi,&&v=\var_i-\vi,
\end{align*}
and choosing $A=O(\d)$ such that $\var_i\geq0$ on $\G_2(D_i)$ and $M_i=\max_{\G_1(B^+_{\rho_1})}\vi$, we deduce from Theorem \ref{t:maxprin} of the Appendix that 
\begin{equation}\label{e:var}
\vi(x)\leq\var_i(x).
\end{equation}
By the Harnack inequality and the assumption that the blow-up is isolated simple, we derive that
\begin{equation}\label{e:mi}
M_i\leq c\vi^{-\l_i}(\yi).
\end{equation}
The estimate \eqref{e:funcbound} of the lemma follows from \eqref{e:var} and \eqref{e:mi}.

To derive \eqref{e:gradbound} from \eqref{e:funcbound}, we argue as follows. For $r_i\leq|\tilde y|\leq\rho_1/2$, we consider
\[
w_i(z)=|\tilde y|^{1-\d}\vi^{\l_i}(\yi)\vi(|\tilde y|z),\quad\textrm{for}\ \frac12\leq|z|\leq2,\ z^3\geq0.
\]
It follows from $(P_i)$ that $w_i$ satisfies 
\begin{equation}\label{e:wi}
\begin{cases}
-\Delta \wi=0,&\textrm{in $\graffe{{\frac12}<|z|<2:\ z^3>0}$,}\\
\displaystyle{\frac{\partial \wi}{\partial \nu}}=H_i(|\tilde y|z)|\tilde y|^{-\l_i}\vi^{\l_i(1-\qi)}(\yi)\wi^{\qi},&\textrm{on $\graffe{{\frac12}<|z|<2:\ z^3=0}$.}
\end{cases}
\end{equation}
In view of \eqref{e:funcbound}, we have $\wi(z)\leq c$ for any $\frac12\leq |z|\leq 2$, $z^3\geq0$. We then derive from~\eqref{e:wi} and gradient elliptic estimates that 
\[
|\n \wi(z)|\leq c,\quad \quad z\in \G_2(B_1^+)
\]
which implies that 
\[
|\n\vi(\tilde y)|\leq c|\tilde y|^{-2+\d}\vi^{-\l_i}(\yi).
\]
This establishes \eqref{e:gradbound}. Estimate \eqref{e:grad2bound} can be derived in a similar way. We omit the details. Lemma \ref{l:bounds} is thus established.
\end{pf}

Later on we will fix $\d$ close to $0$, hence fix $\rho_1$. Our aim is to obtain \eqref{e:funcbound} with $\d=0$ for $r_i\leq d(y,\yi)\leq\rho_1$, which together with Lemma \ref{l:bubble} yields Proposition \ref{p:upperbound}.

Now we state the following Pohozaev type identity, which is basically contained in Li-~Zhu \cite{LZ3}. In the following, we write in some geodesic normal coordinate $x=(x^1,x^2,x^3)$ with $g_{ij}(0)=\d_{ij}$ and $\G_{ij}^k(0)=0$. We use also the notation $\n=(\partial_1,\partial_2,\partial_3)$, $dx\kern0pt=dx^1\wedge dx^2 \wedge dx^3$ and $ds$ to denote the surface area element with respect to the flat metric. 
\begin{lem}\label{l:pohozaev}
For $H\in C^2(\G_1(B_1^+))$ and $a\in C^2(\G_1(B_1^+))$, let $u\in C^2(\ov{B_1^+})$ satisfy, for $q>0$,
\[
\begin{cases}
\displaystyle{-\Delta _gu+\frac 18\, R_g u=0},\quad u>0,&\textrm{in $B_1^+$,}\\
\displaystyle{\frac{\partial u}{\partial \nu}}=Hu^q,&\textrm{on $\G_1(B_1^+)$;}
\end{cases}
\]
then we have, for any $r$ such that $0<r\leq 1$, 
$$\displaylines{\frac1{q+1}\int_{\G_1(B_r^+)}(x'\cdot\n_{x'}H)u^{q+1}\, ds+\tonde{\frac2{q+1}-\frac12}\int_{\G_1(B_r^+)}Hu^{q+1}\,ds\cr
-\frac1{16}\int_{B_r^+}(x\cdot\n R_g)u^2\,dx
-\frac18\int_{B_r^+}R_g\,u^2\,dx
-\frac r{16}\int_{\G_2(B_r^+)}R_g\,u^2\,ds\cr
-\frac r{q+1}\int_{\partial\G_1(B_r^+)}Hu^{q+1}\,ds=\int_{\G_2(B_r^+)}B(r,x,u,\n u)\,ds+A(g,u)\cr
}$$
where
\begin{equation}\label{e:poho1}
B(r,x,u,\n u)=\frac12 \frac{\partial u}{\partial \nu}u+\frac12 r\tonde{\frac{\partial u}{\partial \nu}}^2-\frac12 r|\n_T u|^2,
\end{equation}
$\n_Tu$ denotes the component of $\n u$ which is tangent to $\G_2(B^+_r)$,
\begin{align}\label{e:poho2}
A(g,u)=&\int_{B_r^+}(x^k\partial_ku)(g_{ij}-\d_{ij})\partial_{ij}u\,dx-\int_{B_r^+}(x^l\partial_lu)(g_{ij}-\G^k_{ij}\partial_ku)\,dx\notag\\
&+\frac12\int_{B_r^+}u(g^{ij}-\d^{ij})\partial_{ij}u\,dx-\frac12\int_{B_r^+}u\,g^{ij}\G^k_{ij}\partial_ku\,dx\notag\\
&-\int_{\G_1(B_r^+)}x^i\frac{\partial u}{\partial
  x_i}(g^{ij}-\d^{ij})\frac{\partial u}{\partial
  x_i}\nu_j-\frac{n-2}2\int_{\G_1(B_r^+)}(g^{ij}-\d^{ij})\frac{\partial u}{\partial
  x_i}\nu_j u,
\end{align}
and $\G^k_{ij}$ denotes the Christofell symbol.
\end{lem} 
Regarding the term $A(g,\ui)$, where $\ui$ is a solution of $(P_i)$, we have the following estimate, the proof of which is a direct consequence of Lemma \ref{l:bubble} and Lemma \ref{l:bounds}.
\begin{lem}\label{l:estimateofA}
Let $\{\vi\}_i$ satisfy $(P_i)$, $\yi\to\bar y\in\G_1(B_1^+)$ be an isolated simple blow-up point. Assume $R_i\to+\infty$ and $0<\e_i<e^{-R_i}$ are sequences for which \eqref{e:bubble} and \eqref{e:raggi} hold. Then, for $0<r<\rho_1$, we have
\[
|A(g,\vi)|\leq C_5r\vi^{-2\l_i}(\yi)
\]
where $C_5$ is some constant independent of $i$ and $r$.
\end{lem}
Using Lemma \ref{l:bubble}, Lemma \ref{l:bounds}, Lemma \ref{l:pohozaev}, Lemma \ref{l:estimateofA}, and standard elliptic estimates, we derive the following estimate about the rate of blow-up of the solutions of $(P_i)$.
\begin{lem}\label{l:taui}
Let $\vi$ satisfy $(P_i)$ and $\yi\to\bar y\in\G_1(B_1^+)$ be an isolated simple blow-up point. Assume $R_i\to+\infty$ and $0<\e_i<e^{-R_i}$ are sequences for which \eqref{e:bubble} and \eqref{e:raggi} hold. Then
\[
\tau_i=O\tonde{\vi^{-2\l_i}(y_i)}.
\]
Consequently $\vi^{\tau_i}(\yi)\to 1$ as $i\to\infty$.
\end{lem} 
\begin{lem}\label{l:limsup}
Let $\vi$ satisfy $(P_i)$ and $\yi\to\bar y\in\G_1(B_1^+)$ be an isolated simple blow-up point. Then, for $0<r<\tilde r/2$, we have
\[
\limsup_{i\to+\infty}\max_{y\in\G_2(B^+_r(\yi))}\vi(\yi)\vi(y)\leq C(r).
\]
\end{lem}
\begin{pf}
Due to Lemma \ref{l:harnack}, it is enough to establish the lemma for $r>0$ sufficiently small. Without loss of generality we may take $\bar r=1$. Pick any $y_r\in \G_2(B^+_r)$ and set 
$$\xi_i(y)=\vi^{-1}(y_r)\vi(y),$$
then $\xi_i$ satisfies
\[
\begin{cases}
\displaystyle{-\D_g\xi_i+\frac 18 R_g \xi_i=0},&\textrm{in $B^+_{1/2}(\bar y)$,}\\
\displaystyle{\frac{\partial\xi_i}{\partial\nu}}=H_i\vi^{\qi-1}(y_r)\xi_i^{\qi},&\textrm{on $\G_1(B^+_{1/2}(\bar y))$.}
\end{cases}
\]
It follows from Lemma \ref{l:harnack} that for any compact set $K\subset B^+_{1/2}(\bar y)\setminus\{\bar y\}$, there exists some constant $c(K)$ such that
\[
c(K)^{-1}\leq\xi_i\leq c(K),\quad\textrm{on}\ K.
\]
We also know from \eqref{e:funcbound} that $\vi(y_r)\to0$ as $i\to+\infty$. Then by elliptic theories, we have, after passing to a subsequence, that $\xi_i\to\xi$ in $C^2_{\rm loc}(B^+_{1/2}(\bar y)\setminus\{\bar y\})$, where $\xi$ satisfies
\[
\begin{cases}
\displaystyle{-\Delta_g\xi+\frac 18\, R_g \xi=0},&\textrm{in $B^+_{1/2}(\bar y)$,}\\
\displaystyle{\frac{\partial\xi}{\partial\nu}}=0,&\textrm{on $\G_1(B^+_{1/2})\setminus\{\bar y\}$.}
\end{cases}
\]
From the assumption that $\yi\to\bar y$ is an isolated simple blow-up point of $\{\vi\}_i$, we know that the function $r^{1/2}\bar \xi(r)$ is nonincreasing in the interval $(0,\tilde r)$ and so we deduce that $\xi$ is singular at $\bar y$. So it follows from Corollary \ref{c:corapp} in the Appendix that for $r$ small enough there exists some positive constant $m>0$ independent of $i$ such that for $i$ large we have
\[
-\frac 18\int_{B^+_r} R_g \xi_i =\int_{B^+_r}-\Delta_g\xi_i=-\int_{\G_1(B^+_r)}\frac{\partial\xi_i}{\partial\nu}-\int_{\G_2(B^+_r)}\frac{\partial\xi_i}{\partial\nu}>m-\int_{\G_1(B^+_r)}\frac{\partial\xi_i}{\partial\nu}
\]
which implies that 
\begin{equation}\label{e:lowbound}
-\frac 18\int_{B^+_r} R_g \xi_i+\int_{\G_1(B^+_r)}\frac{\partial\xi_i}{\partial\nu}>m.
\end{equation}
On the other hand
\begin{equation}\label{e:upbound}
\int_{\G_1(B^+_r)}\frac{\partial\xi_i}{\partial\nu}=\int_{\G_1(B^+_r)}H_i\vi^{\qi-1}(y_r )\xi_i^{\qi}\leq\vi^{-1}(y_r)\int_{\G_1(B^+_r)}H_i\vi^{\qi}.
\end{equation}
Using Lemma \ref{l:bubble} and Lemma \ref{l:bounds}, we derive that
\begin{equation}\label{e:finale} 
\int_{\G_1(B^+_r)}H_i\vi^{\qi}\leq c\vi^{-1}(\yi).
\end{equation}
Hence our lemma follows from \eqref{e:lowbound}, \eqref{e:upbound}, and \eqref{e:finale}.
\end{pf}

Now we are able to give the proof of Proposition \ref{p:upperbound}.

\medskip
\begin{pfn}{Proposition \ref{p:upperbound}.} \rm We first establish \eqref{e:upperbound} arguing by contradiction. Suppose the contrary; then, possibly passing to a subsequence still denoted as $\vi$, there exists a sequence $\{\tilde \yi\}_i$ such that $d(\tilde\yi,\yi)\leq\tilde r/2$ and
\begin{equation}\label{e:limite}
\vi(\tilde\yi)\vi(\yi)d(\tilde\yi,\yi)\ \mathop{\longrightarrow}\limits_{i\to+\infty}\ +\infty.
\end{equation}
Set $\tilde r_i=d(\tilde\yi,\yi)$. From Lemma \ref{l:bubble} it is clear that $\tilde r_i\geq r_i=R_i\vi^{1-\qi}(\yi)$. Set
\[
\tilde \vi(x)=\tilde r_i^{\frac1{\qi-1}}\vi(\yi+\tilde r_ix)\quad\textrm{in}\ B_2^{-T_i}, \quad T_i=\tilde r_i^{-1}\yi^3.
\]
Clearly $\tilde v_i$ satisfies
\[
\begin{cases}
\displaystyle{-\D_{g_i}\tilde\vi+\frac 18\, \wtilde R_{g_i}\tilde \vi =0},\quad\vi>0,&\textrm{in}\ B_2^{-T_i},\\
\displaystyle{\frac{\partial\tilde\vi}{\partial\nu}}=\wtilde H_i(x)\tilde\vi^{\qi}(x),&\textrm{on}\ \partial B_2^{-T_i}\cap\{x^3--T_i\},
\end{cases}
\]
where
\begin{align*}
(g_i)_{\a\beta}&= g_{\a\beta}(\tilde r_i x)\,dx^{\a}dx^{\beta},\\
\wtilde R_{g_i}(x)&=\tilde r_i^2 R_{g_i}(\yi+\tilde r_ix),
\end{align*}
and
\[
\wtilde H_i(x)=H_i(\yi+\tilde r_ix).
\]
Lemma \ref{l:limsup} yields that
\[
\max_{x\in \G_2(B^+_{1/2})}\tilde\vi(0)\tilde \vi(x)\leq c
\]
for some positive constant $c$, and so 
\[
\vi(\tilde \yi)\vi(\yi) d(\yi,\yi)\leq c.
\]
This contradicts \eqref{e:limite}. Therefore \eqref{e:upperbound} is established. Take now 
\[
\wi(x)=\vi(\yi)\vi(x).
\]
From $(P_i)$ it is clear that $\wi$ satisfies
\[
\begin{cases}
\displaystyle{-\D_{g}\wi+\frac 18\, R_g w_i=0},&\textrm{in}\ B_3^+,\\
\displaystyle{\frac{\partial\wi}{\partial\nu}}= H_i(x)\vi^{1-\qi}(\yi)\wi^{\qi},&\textrm{on}\ \G_1(B^+_3).
\end{cases}
\]
Estimate \eqref{e:upperbound} implies that $\wi(x)\leq c\,d(x,\yi)^{-1}$. Since $\yi\to \bar y$, $\wi$ is locally bounded in any compact set not containing $\bar y$. Then, up to a subsequence, $\wi\to w$ in $C^2_{\rm loc}(B_{\tilde \rho}(\bar y)\setminus\{ \bar y\})$ for some $w>0$ satisfying
\[
\begin{cases}
\displaystyle{-\D_{g}w+\frac 18\, R_g w=0},&\textrm{in}\ B_{\tilde\rho}^+(\bar y),\\
\displaystyle{\frac{\partial w}{\partial\nu}}= 0,&\textrm{on}\ \G_1(B^+_{\tilde\rho})\setminus\{\bar y\}.
\end{cases}
\]
From Proposition \ref{p:b} of the Appendix, we have that 
\[
w=b\,G(\cdot,\bar y)+E,\quad \textrm{in}\ B^+_{\tilde \rho}\setminus\{0\},
\]
where $b\geq0$, $E$ is a regular function satisfying 
\[
\begin{cases}
\displaystyle{-\D_gE+\frac 18\, R_g E=0},&\textrm{in}\ B^+_{\tilde\rho},\\
\displaystyle{\frac{\partial E}{\partial\nu}}=0,&\textrm{on}\ \G_1(B^+_{\tilde\rho}),
\end{cases}
\]
and $G\in C^2(B^+_{\tilde\rho}\setminus\{\bar y\})$ satisfies
\[
\begin{cases}
-L_gG(\cdot,\bar y)=0,&\textrm{in}\ B^+_{\tilde\rho},\\
\displaystyle{\frac{\partial G_a}{\partial\nu}}=0,&\textrm{on}\ \G_1(B^+_{\tilde\rho})\setminus\{\bar y\},
\end{cases}
\]
and $\lim_{y\to\bar y}d(y,\bar y)G(y,\bar y)$ is a constant. Moreover $w$ is singular at $\bar y$. Indeed from the definition of isolated simple blow-up we know that the function $r^{1/2}\bar w(r)$ is nonincreasing in the interval $(0,\tilde r)$, which implies that $w$ is singular at the origin and hence $b>0$. The proof of Proposition \ref{p:upperbound} is thereby complete.
\end{pfn}

Using Proposition \ref{p:upperbound}, one can strengthen the results of Lemmas \ref{l:bounds} and \ref{l:estimateofA} just using \eqref{e:upperbound} instead of \eqref{e:funcbound}, thus obtaining the following corollary.
\begin{cor}\label{c:coruppbound}
Let $\{\vi\}_i$ satisfy $(P_i)$, $\yi\to\bar y\in\G_1(B_1^+)$ be an isolated simple blow-up point. Assume $R_i\to+\infty$ and $0<\e_i<e^{-R_i}$ are sequences for which \eqref{e:bubble} and \eqref{e:raggi} hold. Then there exists $\rho_1\in(0,\tilde r)$ such that 
\begin{equation}\label{e:gradbound2}
|\n_g\vi(y)|\leq C_4\vi^{-1}(\yi)d(y,\yi)^{-2},\quad\textrm{for all}\ r_i\leq d(y,\yi)\leq\rho_1,
\end{equation}
and
\begin{equation}\label{e:grad2bound2}
|\n^2_g\vi(y)|\leq C_4\vi^{-1}(\yi)d(y,\yi)^{-3},\quad\textrm{for all}\ r_i\leq d(y,\yi)\leq\rho_1,
\end{equation}
where $r_i=R_i\vi^{1-\qi}(\yi)$ and $C_4$ is some positive constant independent of $i$. Moreover
\[
|A(g,\vi)|\leq C_5 r\vi^{-2}(\yi),
\]
for some positive constant $C_5$ independent of $i$.
\end{cor}
Let us prove an upper bound estimate for $\n_g H_i(\yi)$.
\begin{lem}\label{l:H}
Let $\vi$ satisfy $(P_i)$ and $\yi\to\bar y\in \G_1(B_1^+)$ be an isolated simple blow-up point. Then
\[
\n_g H_i(\yi)=O(\vi^{-2}(\yi)).
\]
\end{lem}
\begin{pf}
Let $x=(x^1,x^2,x^3)$ be some geodesic normal coordinates centered at $\yi$ and $\eta$ some smooth cut-off function such that 
\[
\eta(x)=
\begin{cases}
1,&\textrm{if}\ x\in\ov{B^+_{1/4}},\\
0,&\textrm{if}\ x\not\in\ov{B^+_{1/2}},
\end{cases}
\]
and $0\leq\eta\leq1$. Multiply $(P_i)$ by $\eta\frac{\partial\vi}{\partial x_1}$ and integrate by parts over $B_1^+$, thus obtaining 
\begin{equation}\label{e:H1}
0=\int_{B_1^+}\n_g\vi\cdot\n\Big(\eta\frac{\partial\vi}{\partial
  x_1}\Big)\,dV+\frac 18 \int_{B_1^+}R_g \vi \eta \frac{\partial\vi}{\partial
  x_1}-\int_{\G_1(B^+_{1/2})}\frac{\partial\vi}{\partial\nu}\eta\frac{\partial\vi}{\partial x_1}\,d\s.
\end{equation}
From $(P_i)$, \eqref{e:upperbound}, and \eqref{e:bubble} we have that
\begin{eqnarray}\label{e:H2}
&&\hskip-1.5truecm\int_{\G_1(B_{1/2}^+)}\frac{\partial\vi}{\partial\nu}\eta\frac{\partial\vi}{\partial x_1}\,d\,\s+\frac 18 \int_{B_1^+}R_g \vi \eta \frac{\partial\vi}{\partial
  x_1}\notag\\
&=&\int_{\G_1(B_{1/2}^+)}H_i\vi^{\qi}\eta\frac{\partial\vi}{\partial x_1}\,d\s+O(\vi^{-2}(\yi))\notag\\
&=&-\frac{1}{\qi+1}\frac{\partial H_i}{\partial x_1}(\yi)\int_{\G_1(B_{1/2}^+)}\eta\vi^{\qi+1}\,d\s+O\bigg(\int_{\G_1(B_{1/2}^+)}|x'|\vi^{\qi+1}\bigg)+O(\vi^{-2}(\yi))\notag\\
&=&-\frac{1}{\qi+1}\frac{\partial H_i}{\partial x_1}(\yi)\int_{\G_1(B_{1/2}^+)}\eta\vi^{\qi+1}\,d\s+O(\vi^{-2}(\yi)).
\end{eqnarray}
On the other hand, from \eqref{e:gradbound2} it follows that
\begin{eqnarray}\label{e:H3}
&&\hskip-1.5truecm\int_{\G_1(B_1^+)}\n_g\vi\cdot\n_g\Big(\eta\frac{\partial\vi}{\partial x_1}\Big)\,d\s\notag\\
&=&\int_{B_1^+}(\n_g\vi\cdot\n_g\eta)\frac{\partial\vi}{\partial x_1}\,dV+\int_{B_1^+}\n_g\vi\cdot\eta\n_g\Big(\frac{\partial\vi}{\partial x_1}\Big)\,dV\notag\\
&=&-\frac12\int_{B^+_{1/2}\setminus B^+_{1/4}}\frac{\partial\eta}{\partial x_1}|\n_g\vi|^2\,dV+O(\vi^{-2}(\yi))=O(\vi^{-2}(\yi)).
\end{eqnarray}
Putting all together \eqref{e:H1}, \eqref{e:H2}, and \eqref{e:H3}, we find
\[
\frac{\partial H_i}{\partial x_1}(\yi)=O(\vi^{-2}(\yi)).
\]
Repeating the same argument for the derivatives with respect to $x_2$ and $x_3$, we come to the required estimate.
\end{pf}
\begin{cor}\label{c:H}
Under the same assumptions of Lemma \ref{l:H}, one has that
\[
\int_{\G_1(B^+_r)}x'\cdot\n_{x'}H_i\vi^{\qi+1}\,d\s=O(\vi^{-4}(\yi)).
\]
\end{cor}
\begin{pf}
We have that
$$
\displaylines{\int_{\G_1(B^+_r)}x'\cdot\n_{x'}H_i\vi^{\qi+1}\,d\s=\int_{\G_1(B^+_r)}\n_{x'}H_i(\yi)\cdot (x'-\yi)\vi^{\qi+1}\,d\s\cr
+O\bigg(\int_{\G_1(B^+_r)}|x'|^2\vi^{\qi+1}\,d\s\bigg).}
$$
Since, using Proposition \ref{p:upperbound} and Lemma \ref{l:bubble}, $\int_{\G_1(B^+_r)}(x'-\yi)\vi^{\qi+1}\,d\s=O(\vi^{-2}(\yi))$, from the previous lemma, Corollary \ref{c:coruppbound}, and \eqref{e:bubble}, we reach the conclusion.
\end{pf}
 
\begin{pro}\label{p:A}
Let $\vi$ satisfy $(P_i)$, $\yi\to\bar y$ be an isolated simple blow-up point with, for some $\tilde\rho>0$,
\[
\vi(\yi)\vi\ \mathop{\longrightarrow}\limits_{i\to+\infty}\ h,\quad\textrm{in}\ C^2_{\rm loc}(B^+_{\tilde\rho}(\bar y)\setminus\{\bar y\}).
\]
Assume, for some $\beta>0$, that in some geodesic normal coordinate system $x=(x^1,x^2,x^3)$
\[
h(x)=\frac\beta{|x|}+A+o(1), \quad\textrm{as}\ |x|\to0.
\]
Then $A\leq0$.
\end{pro}
\begin{pf}
For $r>0$ small, the Pohozaev type identity of Lemma \ref{l:pohozaev} yields
\begin{align}\label{e:pohoz}
&\frac1{\qi+1}\int_{\G_1(B^+_r)}(x'\cdot\n_{x'}H_i)\vi^{\qi+1}\,ds+\Big(\frac2{\qi+1}-\frac12\Big)\int_{\G_1(B^+_r)}H_i\vi^{\qi+1}\,ds\notag\\
&-\frac1{16}\int_{B_r^+}(x\cdot\n R_g)\vi^2\,dx
-\frac18\int_{B_r^+}R_g\,\vi^2\,dx
-\frac r{16}\int_{\G_2(B_r^+)}R_g\,\vi^2\,ds\notag\\
&\quad-\frac r{\qi+1}\int_{\partial\G_1(B^+_r)}H_i\vi^{\qi+1}=\int_{\G_2(B^+_r)}B(r,x,\vi,\n\vi)\,ds+A(g,\vi)
\end{align}     
where $B$ and $A(g,\vi)$ are defined in \eqref{e:poho1} and \eqref{e:poho2} respectively. Multiply \eqref{e:pohoz} by $\vi^2(\yi)$ and let $i\to\infty$. Using Corollary \ref{c:coruppbound}, Lemma \ref{l:bubble}, and Corollary \ref{c:H}, one has that 
\begin{equation}\label{e:pos}
\lim_{r\to 0^+}\int_{\G_2(B^+_r)}B(r,x,h,\n h)=\lim_{r\to 0^+}\limsup_{i\to \infty}\vi^2(\yi)\int_{\G_2(B^+_r)}B(r,x,\vi,\n \vi)\geq0.
\end{equation}
On the other hand, a direct calculation yields
\begin{equation}\label{e:ug}
\lim_{r\to0^+}\int_{\G_2(B^+_r)}B(r,x,h,\n h)=-c\,A
\end{equation}
for some $c>0$. The conclusion follows from \eqref{e:pos} and \eqref{e:ug}.
\end{pf}

Now we can prove that an isolated blow-up point is in fact an isolated simple blow-up point.
\begin{pro}\label{p:is}
Let $\vi$ satisfy $(P_i)$ and $\yi\to\bar y$ be an isolated blow-up point. Then $\bar y$ must be an isolated simple blow-up point.
\end{pro}
\begin{pf}
The proof is basically the same as that of Proposition 2.11 of \cite{VA}. For the reader's convenience, we include the proof here. From Lemma \ref{l:bubble}, it follows that 
\begin{equation}\label{e:is1}
\bar\wi'(r)<0\quad\textrm{for every}\ C_2\vi^{1-\qi}(\yi)\leq r\leq r_i.
\end{equation}
Suppose that the blow-up is not simple; then there exist some sequences $\tilde r_i\to 0^+$, $\tilde c_i\to~+\infty$ such that $\tilde c_i\vi^{1-\qi}(\yi)\leq\tilde r_i$ and, after passing to a subsequence, 
\begin{equation}\label{e:is2}
\bar w_i'(\tilde r_i)\geq 0.
\end{equation}
From \eqref{e:is1} and \eqref{e:is2} it is clear that $\tilde r_i\geq r_i$ and $\bar w_i$ has at least one critical point in the interval $[r_i,\tilde r_i]$. Let $\mu_i$ be the smallest critical point of $\bar w_i$ in this interval. We have that 
\[
\tilde r_i\geq\mu_i\geq r_i\quad\textrm{and}\quad\lim_{i\to\infty}\mu_i=0.
\]
Let $g_i=(g_i)_{\a\beta}\,dx^{\a}dx^\beta=g_{\a\beta}(\mu_i x)\,dx^{\a}dx^\beta$ be the scaled metric and
\[
\xi_i(x)=\mu_i^{\frac1{\qi-1}}\vi(\yi+\mu_i x).
\]
Then $\xi_i$ satisfies
\begin{equation}\label{e:lu}
\begin{cases}
\displaystyle{-\Delta_{g_i}\xi_i+\frac 18\, R_{g_i}\xi_i=0},&\textrm{in}\ B^{-T_i}_{1/\mu_i},\\
\displaystyle{\frac{\partial \xi_i}{\partial \nu}}=\wtilde H_i(x)\xi^{\qi}\xi_i^{\qi},&\mbox{on}\ \partial B^{-T_i}_{1/\mu_i}\cap\{x^3=-T_i\},\\
\lim_{i\to\infty}\xi_i(0)=\infty&\mbox{and}\ 0\ \mbox{is a local maximum point of}\ \xi_i,\\
r^\frac{1}{\qi-1}\bar\xi_i(r)&\hbox{has negative derivative in}\quad c\,\xi_i(0)^{1-\qi}<r<1,\\
\left.\frac{d}{dr}\left(r^\frac{1}{\qi-1}\bar\xi_i(r)\right)\right|_{r=1}=0,
\end{cases}
\end{equation}
where $T_i=\mu_i^{-1}\yi^3$, $\tilde a_i(x)=\mu_ia_i(\yi+\mu_i x)$, and $\wtilde H_i(x)=H_i(\yi+\mu_i x)$. Arguing as in the proof of Lemma \ref{l:bubble}, we can easily prove that $T_i\to 0$. Since $0$ is an isolated simple blow-up point, by Proposition \ref{p:upperbound} and Lemma \ref{l:harnack}, we have that, for some $\b>0$, 
\begin{equation}\label{e:conv}
\xi_i(0)\xi_i\ \mathop{\longrightarrow}\limits_{i\to+\infty}\ h=\b\, |x|^{-1}+E\quad\textrm{in}\ C^2_{\rm loc}(\R^3_+\setminus\{0\})
\end{equation}
with $E$ satisfying
\[
\begin{cases}
-\D E=0,&\textrm{in}\ \R^3_+,\\
\displaystyle{\frac{\partial E}{\partial\nu}}=0,&\textrm{on}\ \partial\R^3_+.
\end{cases}
\]
By the Maximum Principle we have that $E\geq0$. Reflecting $E$ to be defined on all $\R^3$ and thus using the Liouville Theorem, we deduce that $E$ is a constant. Using the last equality in \eqref{e:lu} and \eqref{e:conv}, we deduce that $E\equiv b$. Therefore, $h(x)=b(G_a(x,\bar y)+1)$ and this fact contradicts Proposition \ref{p:A}.
\end{pf}

\section{Ruling out bubble accumulations}
\mbox{}

Now we can proceed as in \cite{VA} to obtain the following results which rule out the possible accumulations of bubbles, and this implies that only isolated blow-up points may occur to blowing-up sequences of solutions.
\begin{pro}\label{p:acc}
 Let $(M,g)$ be a smooth compact three dimensional Riemannian
 ma\-ni\-fold with umbilic boundary. For any $R\geq1$, $0<\e<1$, there exist positive constants $\delta_0$, $c_0$, and $c_1$ depending only on $M$, $g$, $\|H\|_{C^2(\partial M)}$, $\inf_{y\in\partial M}H(y)$, $R$, and~$\e$, such that for all $u$ in
$$\bigcup_{3-\delta_0\leq q\leq 3}{\cal M}_{H,q}$$
with $\max_Mu\geq c_0$, there exists ${\cal S}=\{p_1,\dots,p_N\}\subset \partial M$ with $N\geq 1$ such that 
\begin{itemize}
\item[\rm(i)] each $p_i$ is a local maximum point of $u$ in $M$ and
\[
\overline{B_{\bar r_i}(p_i)}\cap\overline{B_{\bar r_j}(p_j)}=\emptyset,\quad\mbox{for}\ i\not=j,
\]
where $\bar r_i=Ru^{1-q}(p_i)$ and $\overline{B_{\bar r_i}(p_i)}$ denotes the geodesic ball in $(M,g)$ of radius $\bar r_i$ and centered at $p_i$;
\item[\rm(ii)] \[
\left\|u^{-1}(p_i)u(\exp_{p_i}(yu^{1-q}(p_i))-\left(\frac{1}{(1+hx^3)^2+h^2|x'|^2}\right)^{1/2}\right\|_{C^2(B_{2R}^M(0))}<\e
\]
where 

\[
B_{2R}^M(0)=\{y\in T_{p_i}M:\ |y|\leq 2R,\ u^{1-q}(p_i)y\in\exp^{-1}_{p_i}(B_\delta(p_i))\},
\]
 $y=(y',y^n)\in\R^n$, and $h>0$;
\item[\rm(iii)] $d^\frac{1}{q-1}(p_j,p_i)u(p_j)\geq c_0$, for $j>i$, while $d(p,{\cal S})^{\frac{1}{q-1}}u(p)\leq c_1$, $\forall \, p\in M$, where $d(\cdot,\cdot)$ denotes the distance function in metric $g$.
\end{itemize}
\end{pro}
\begin{pro}\label{p:noaccumulation}
Let $(M,g)$ be a smooth compact three dimensional Riemannian manifold
with umbilic boundary. For suitably large $R$ and small $\e>0$, there exist $\delta_1$ and $d$ depending only on $M$, $g$, $\|a\|_{C^2(\partial M)}$, $\|H\|_{C^2(\partial M)}$, $\inf_{y\in\partial M}H(y)$, $R$, and $\e$, such that for all $u$ in
$$\bigcup_{3-\delta_1\leq q\leq3}{\cal M}_{a,H,q}$$
with $\max_Mu\geq c_0$, we have
$$\min\{d(p_i,p_j):\ i\not=j,\ 1\leq i,j\leq N\}\geq d$$
where $c_0,p_1,\dots,p_N$ are given by Proposition \ref{p:acc}.
\end{pro}
The previous two propositions imply that any blow-up point is in fact an isolated blow-up point. Thanks to Proposition \ref{p:is}, any blow-up point is in fact an isolated simple blow-up~point.

\medskip
\begin{pfn}{Theorem \ref{t:sth1}.} Arguing by contradiction, suppose that there exist some sequences $\qi\to q\in]1,3]$, $\ui\in{\cal M}_{H_i,\qi}$ such that $\|\ui\|_{H^1(M)}\to+\infty$ as $i\to\infty$, which, in view of standard elliptic estimates, implies that $\max_{M}\ui\to+\infty$. 

From Hu \cite{hu} (see also \cite{lizhang}), we know that $q=3$. By Proposition \ref{p:noaccumulation}, we have that for some small $\e>0$, large $R>0$, and some $N\geq1$ there exist $\yi^{(1)},\dots,\yi^{(N)}\in\partial M$ such that (i-iii) of Proposition \ref{p:acc} hold. $\{\yi^{(1)}\}_i,\dots,\{\yi^{(N)}\}_i$ are isolated blow-up points and hence, by Proposition \ref{p:is}, isolated simple blow-up points. From \eqref{e:bubble} and Proposition \ref{p:upperbound}, we have that $\{\|\ui\|_{H^1(M)}\}_i$ is bounded, thus finding a contradiction. Theorem \ref{t:sth1} is thereby established.
\end{pfn}

\section{Compactness of the solutions}
\mbox{}

Before proving Theorem \ref{t:pos}, we state the following result about the compactness of solutions of $(P_{H,q})$ when $q$ stays strictly below the critical exponent. The proof is basically the same as the proof of Theorem 3.1 of \cite{VA}.
\begin{thm}\label{t:below}
Let $(M,g)$ be a smooth compact three dimensional Riemannian manifold
with umbilic boundary. Then for any $\delta_1>0$ there exists a constant $C>0$ depending only on $M$, $g$, $\delta_1$, $\|H\|_{C^2(\partial M)}$, and the positive lower bound of $H$ on $\partial M$ such that for all $u\in\bigcup_{1+\delta_1\leq q\leq3-\delta_1}{\cal M}_{H,q}$ we have
\[
\frac{1}{C}\leq u(x)\leq C,\quad\forall\, x\in M;\quad\|u\|_{C^2(M)}\leq C.
\]
\end{thm}
Now we prove Theorem \ref{t:pos}.

\medskip
\begin{pfn}{Theorem \ref{t:pos}.} Due to elliptic estimates and
  Lemma \ref{l:harnack}, we have to prove just the $L^\infty$ bound,
  i.e. $u\leq C$. Suppose the contrary; then there exist a
  sequence $\qi\to q\in]1,3]$ with 
\[
\ui\in{\cal M}_{H,\qi},\quad\textrm{and}\ \max_M\ui\to+\infty,
\]
where $\bar c$ is some positive constant independent of $i$. From Theorem \ref{t:below}, we have that $q$ must be $3$. It follows from Proposition \ref{p:is} and  Proposition \ref{p:noaccumulation} that, after passing to a subsequence, $\{\ui\}_i$ has $N$ ($1~\leq~N<~\infty$) isolated simple blow-up points denoted by $y^{(1)},\dots,y^{(N)}$. Let $y_i^{(\ell)}$ denotes the local maximum points as in Definition \ref{d:isolated}. It follows from Proposition \ref{p:upperbound} that 
\[
\ui(\yi^{(1)})\ui\ \mathop{\longrightarrow}\limits_{i\to+\infty}\ h(y)=\sum_{j=1}^Nb_jG(y,y^{(j)})+E(y)\quad\textrm{in}\ C^2_{\rm loc}\big(M\setminus \{y^{(1)},\dots,y^{(N)}\}\big),
\]
where $b_j>0$ and $E\in C^2(M)$ satisfies
\begin{equation}\label{e:eqofE}
\begin{cases}
-L_gE=0,&\textrm{in}\ M,\\
\displaystyle{\frac{\partial E}{\partial\nu}}=0,&\textrm{on}\ \partial M.
\end{cases}
\end{equation}
Since the manifold is of positive type we have that $E\equiv0$. Therefore,
\[
\ui(\yi^{(1)})\ui\ \mathop{\longrightarrow}\limits_{i\to+\infty}\ h(y)=\sum_{j=1}^Nb_jG_a(y,y^{(j)})\quad\textrm{in}\ C^2_{\rm loc}\big(M\setminus \{y^{(1)},\dots,y^{(N)}\}\big).
\]
Let $x=(x^1,x^2,x^3)$ be some geodesic normal coordinate system
centered at $\yi^{(1)}$. From Lemma \ref{l:G} of the Appendix, the
Positive Mass Theorem, and the assumption that the manifold is not
conformally equivalent to the standard ball, we derive that there exists a positive constant $A$ such that 
\[
h(x)=h(\exp_{\yi^{(1)}}(x))=c|x|^{-1}+A_i+O(|x|^{-\a})\quad\textrm{for $|x|$ close to $0$}
\]
and $A_i\geq A>0$. This contradicts the result of Proposition \ref{p:A}. The compactness part of Theorem \ref{t:pos} is proved.

Since we have compactness, we can proceed as in section 4 of \cite{VA} to prove that the total degree of the solutions is $-1$. Theorem \ref{t:pos} is established.
\end{pfn}

\appendix
\section*{Appendix}
\setcounter{section}{1}
\setcounter{equation}{0}
\setcounter{thm}{0}

\mbox{}

In this Appendix, we recall some well known results and provide some description of singular behaviour of positive solutions to some boundary value elliptic equations in punctured half balls.\\\indent
For $n\geq 3$ let $B_r^+$ denote the set $\{x=(x',x^n)\in\R^n=\R^{n-1}\times\R:\ |x|<r\ \mbox{and}\ x^n>~0\}$ and set $\G_1(B_r^+):=\partial B^+_r\cap\partial\R^n_+$, $\G_2(B_r^+):=\partial B^+_r\cap\R^n_+$. Throughout this section, let $g= g_{ij}\,dx^i dx^j$ denote some smooth Riemannian metric in $B_1^+$ and $a\in C^1(\G_1(B_1^+))$.

First of all we recall the following Maximum Principle; for the proof see \cite{HL1}.
\begin{thm}\label{t:maxprin}
Let $\Omega$ be a bounded domain in $\R^n$ and let $\partial \Omega=\Gamma\cup\Sigma$, $V\in L^{\infty}(\Omega)$, and $h\in L^\infty(\Sigma)$ such that there exists some $\psi\in C^2(\Omega)\cap C^1(\ov{\Omega})$, $\psi>0$ in $\ov{\Omega}$ satisfying
\[
\begin{cases}
\Delta_g\psi+V\psi\leq 0,&\mbox{in}\ \Omega,\\
\displaystyle{\frac{\partial\psi}{\partial\nu}}\geq h\psi,&\mbox{on}\ \Sigma.
\end{cases}
\]
If $v\in C^2(\Omega)\cap C^1(\ov\Omega)$ satisfies
\[
\begin{cases}
\Delta_g v+Vv\leq 0,&\mbox{in}\ \Omega,\\
\displaystyle{\frac{\partial v}{\partial\nu}}\geq hv,&\mbox{on}\ \Sigma,\\
v\geq 0,&\mbox{on}\ \Gamma,
\end{cases}
\]
then $v\geq 0$ in $\ov\Omega$.
\end{thm}

We state now the following Maximum Principle which holds for the operator $T$ defined by
\[
Tu=v\quad\textrm{if and only if}\quad
\begin{cases}
L_gu=0,&\textrm{in}\ \mathring M,\\
\displaystyle{\frac{\partial u}{\partial\nu}}=v,&\textrm{on}\ \partial M.
\end{cases}
\]
\begin{pro}\label{p:mp}
{\bf \cite{E2}} Let $(M,g)$ be a Riemannian manifold with boundary of
positive type. Then for any $u\in C^2(\mathring M)\cap C^1(M)$ satisfying
\[
\begin{cases}
L_gu\geq 0,&\textrm{in}\ \mathring M,\\
\displaystyle{\frac{\partial u}{\partial\nu}}\leq 0,&\textrm{on}\ \partial M,
\end{cases}
\]
we have $u\leq 0$ in $M$.
\end{pro}
\begin{pf} Let $u^+(x)=\max\{0,u(x)\}$. Then 
\begin{equation*}
0\leq\int_M (L_g u)u^+\,dV-\int_{\partial M}\frac{\partial u}{\partial\nu}u^+\,d\s
=-\int_M|\n_g u^+|^2\,dV-\frac 18\int_MR_g^2|u^+|^2\,dV.
\end{equation*}
Sine $M$ is of positive type $\int |\n _g u|^2+\frac 18 \int R_g u^2$
is an equivalent norm hence $u^+\equiv 0$.
\end{pf}

\noindent We now recall the following Louville type Theorem by Li and Zhu \cite{LZ1}.
\begin{thm}\label{t:liouville}
If $v$ is a solution of
\[
\begin{cases}
-\Delta v=0,&\textrm{in}\ \R^n_+,\\
\displaystyle{\frac{\partial v}{\partial x^n}}=c v^{\frac{n}{n-2}},&\textrm{on}\ \partial\R^n_+,
\end{cases}
\]
and $c$ is a negative constant, then either $v\equiv0$ or $v$ is of the form
\[
v(x',x^n)=\left[\frac{\e}{(x_0^n+x^n)^2+|x'-x_0'|^2}\right]^\frac{n-2}{2},\quad x'\in\R^{n-1},\ x^n\in\R,
\]
where $x_0^n=-\frac{(n-2)\e}{c}$, for some $\e>0$, and $x_0'\in\R^{n-1}$.
\end{thm}
\begin{lem}\label{l:reg}
Suppose that $u\in C^2(B_1^+\setminus\{0\})$ is a solution of
\begin{equation}\label{e:reg}
\begin{cases}
-L_gu=0,&\mbox{on $B_1^+$,}\\
\displaystyle{\frac{\partial u}{\partial\nu}}=0
,&\mbox{on $\G_1(B_1^+\setminus\{0\})$,}
\end{cases}
\end{equation}
and $u(x)=o\tonde{|x|^{2-n}}$ as $|x|\to0$. Then $u\in C^{2,\a}(B^+_{1/2})$ for any $0<\a<1$.
\end{lem}
\begin{pf}
We make a reflection cross $\Gamma_1(B_1^+)$ to extend $u$ as a
solution of $-L_g u=0$ on $B_1^+$, then we use \cite{GS} to conclude
that $0$ is a removable singularity, then the result follows from
standard elliptic regularity.
\end{pf}

\begin{lem} \label{l:G}
There exists some constant $\d_0>0$ depending only on $n$,
$\|g_{ij}\|_{C^2(B_1^+)}$ and $\|H\|_{L^{\infty}(B_1^+)}$ such that
for all $0<\d<\d_0$ there exists some function $G$ satisfying 
\begin{equation}\label{e:G}
\begin{cases}
-L_gG=0,&\textrm{in $B^+_{\d}$,}\\[5pt]
\displaystyle{\frac{\partial G}{\partial\nu}}=0,&\textrm{on $\G_1(B^+_{\d})\setminus\{0\}$,}\\
\lim_{|x|\to 0}|x|^{-1}G(x)=1
\end{cases}
\end{equation}
such that, for some $A$ constant and some $\a\in (0,1)$, $\forall\, x\in B_{\d}^+$ 
\[
G(x)=|x|^{-1}+A+O(|x|^{\a}).
\]
\end{lem}
\begin{pf}
Reflecting across $\Gamma_1(B_{\d}^+)$, the lemma is reduced to
Proposition B.1
 in \cite{LZ3}
\end{pf}

Making again reflection across $\Gamma_1(B_{1}^+)$, we derive from Lemma 9.3
in \cite{LZ3}, the following
\begin{lem}\label{l:alfa}
Assume that $u\in C^2(B^+_1\setminus\{0\})$ satisfies
\[
\begin{cases}
-L_gu=0,&\textrm{in $B^+_1$,}\\
\displaystyle{\frac{\partial u}{\partial\nu}}=0,&\textrm{on $\G_1(B^+_1)\setminus\{0\}$,}
\end{cases}
\]
then
\[
\a=\limsup_{r\to 0^+}\max_{x\in\G_2(B_r^+)}u(x)|x|^{n-2}<+\infty.
\]
\end{lem}

\begin{pro}\label{p:b}
Suppose that $u\in C^2(B_1^+\setminus\{0\})$ satisfies 
\begin{equation}\label{e:equ}
\begin{cases}
-L_gu=0,&\textrm{in $B_1^+$,}\\
\displaystyle{\frac{\partial u}{\partial \nu}}=0,&\textrm{on $\G_1(B^+_1)\setminus\{0\}$.}
\end{cases}
\end{equation}
Then there exists some constant $b\geq0$ such that
\[
u(x)=b\,G(x)+E(x)\quad\textrm{in}\ B^+_{1/2}\setminus\{0\},
\]
where $G$ is defined in Lemma \ref{l:G}, and $E\in C^2(B_1^+)$ satisfies
\begin{equation}\label{e:eqe}
\begin{cases}
-L_gE=0,&\textrm{in $B_{1/2}^+$,}\\
\displaystyle{\frac{\partial E}{\partial \nu}}=0,&\textrm{on $\G_1(B^+_{1/2})$.}
\end{cases}
\end{equation}
\end{pro}
\begin{pf}
Set 
\begin{equation}\label{e:b}
b=b(u)=\sup\big\{\l\geq0:\ \l G\leq u\ \textrm{in}\ \ov{B^+_{\d_0}}\setminus\{0\}\big\}.
\end{equation}
By the previous lemma we know that 
\[
0\leq b\leq \a<+\infty.
\]
Two cases may occur.
\begin{description}
\item[Case 1:] $b=0$.\\ In this case we claim that 
\[
\foral \e>0\quad \exists\,r_\e\in(0,\d_0):\ \min_{x\in\G_2(B^+_r)}\{u(x)-\e G(x)\}\leq0,\ \textrm{for any}\  0< r<r_\e.
\]
To prove the claim argue by contradiction. Suppose that the claim is false. Then there exist $\e_0>0$ and a sequence $r_j\to 0^+$ such that
\[
\min_{|x|=r_j}\{u(x)-\e_0 G(x)\}>0\quad\textrm{and}\quad u(x)-\e_0 G(x)>0\ \textrm{on}\ \G_2(B^+_{\d_0}).
\]
Let us prove that $\e_0\leq b$, which gives a contradiction. To do this, we want to prove that 
\[
u-\e_0 G\geq 0\quad\textrm{in}\ B^+_{\d_0}\setminus\{0\}.
\]
Note that $u-\e_0 G$ satisfies
\[
\begin{cases}
\displaystyle{-\D_g(u-\e_0G)+\frac18\, R_g(u-\e_0 G)=0},&\textrm{in}\ B^+_{\d_0},\\
\displaystyle{\frac{\partial }{\partial \nu}}(u-\e_0 G)=0,&\textrm{on}\ \G_1(B^+_{\d_0})\setminus\{0\}.
\end{cases}
\]
Apply Theorem \ref{t:maxprin} with
\[
v=u-\e_0G,\quad \Sigma=\G_1(B^+_{\d_0})\setminus\G_1(B^+_{r_j}),\quad \G=\G_2(B^+_{r_j})\cup \G_2(B^+_{\d_0}),
\]
thus getting $u-\e_0G\geq0$ in the annulus $B^+_{\d_0}\setminus B^+_{r_j}$ for any $j$, and consequently in $B^+_{\d_0}\setminus \{0\}$. Therefore $\e_0\leq b$ and so there is a contradiction. This proves the claim. Hence, for any $\e>0$ and $0<r<r_\e$ there exists $x_\e\in\G_2(B^+_r)$ such that
\[
u(x_\e)\leq \e G(x_\e).
\]
By the Harnack inequality of Lemma \ref{l:harnack} we have that 
\[
\max_{|x|=r}u(x)\leq c\,u(x_\e)\leq c\,\e G(x_\e).
\]
Since $G(x)\sim|x|^{2-n}$ for $|x|$ small, we conclude that
\[
u(x)=o\big(|x|^{2-n}\big)\quad \textrm{for}\ |x|\sim0.
\]
Therefore from Lemma \ref{l:reg} we obtain that $u$ is regular. Setting $E(x)=u(x)$, the conclusion in this case follows.
\item[Case 2:] $b>0$.\\ We consider $v(x)=u(x)-b\,G(x)$ in $\ov{B^+_{\d_0}}\setminus\{0\}$. By definition of $b$, it is clear that $v\geq0$ in $\ov{B^+_{\d_0}}\setminus\{0\}$. Moreover $v$ satisfies
\[
\begin{cases}
\displaystyle{-\D_gv+\frac18\, R_g v=0},&\textrm{in $B_{\d_0}^+$,}\\
\displaystyle{\frac{\partial v}{\partial \nu}}=0,&\textrm{on $\G_1(B^+_{\d_0})\setminus\{0\}$,}
\end{cases}
\]
so that from the Maximum Principle we know that either $v\equiv 0$ or $v>0$ in $B_{\d_0}^+\setminus\{0\}$. If $v\equiv0$, take $E\equiv 0$ and we are done. Otherwise $v>0$ and satisfies the same equation as $u$. Set 
\[
\tilde b=b(v)=\sup\big\{\l\geq0:\ \l G\leq v\ \textrm{in}\ \ov{B^+_{\d_0}}\setminus\{0\}\big\}.
\]
If $\l\geq0$ and $\l G\leq v$ in $B^+_{\d_0}\setminus \{0\}$, then $\l G\leq u-b\,G$ with $b>0$, i.e. $(\l+b)G\leq u$. By definition of $b$, this implies that $\l+b\leq b$, i.e. $\l\leq 0$ and so $\l=0$. Therefore $\tilde b=0$. Arguing as in Case 1, we can prove that $v(x)=o\big(|x|^{2-n}\big)$ for $|x|\sim 0$. Lemma \ref{l:reg} ensures that $v$ is regular, so that choosing $E(x)=v(x)$ we are done.
\end{description}
The proof of Proposition \ref{p:b} is thereby complete.
\end{pf}

\noindent Finally we prove the following corollary.
\begin{cor}\label{c:corapp}
Let $u$ be a solution of \eqref{e:equ} which is singular at $0$. Then
\[
\lim_{r\to0^+}\int_{\G_2(B^+_r)}\frac{\partial u}{\partial \nu}\,d\s=b\cdot\lim_{r\to0^+}\int_{\G_2(B^+_r)}\frac{\partial G}{\partial \nu}\,d\s=-\frac{n-2}2\,b\,|{\mathbb S}^{n-1}|,
\]
where ${\mathbb S}^{n-1}$ denotes the standard $n$-dimensional sphere
and $b>0$ is given by Proposition~\ref{p:b}.
\end{cor}
\begin{pf}
From the previous proposition, we know that 
\[
u(x)=b\,G(x)+E(x),\quad\textrm{in}\ B^+_{1/2}(0)\setminus\{0\},\quad b\geq0.
\]
Since $u$ is singular at $0$, $b$ must be strictly positive. From \eqref{e:eqe}, we have
\[
0=-\int_{B^+_r}\D_gE\,dV-\int_{\G_2(B_r^+)}\frac{\partial E}{\partial
  \nu}\,d\s+\frac 18\int_{B^+_r}R_g E .
\]
Hence, since $E$ is regular, we obtain 
\[
\int_{\G_2(B_r^+)}\frac{\partial E}{\partial \nu}\,d\s=\frac 18\int_{B_r^+}R_g(x)E(x)\,d\s\ \mathop{\longrightarrow}\limits_{r\to0^+}\ 0
\]
and so 
\[
\lim_{r\to0^+}\int_{\G_2(B^+_r)}\frac{\partial u}{\partial\nu}\,d\s=\lim_{r\to0^+}b\int_{\G_2(B^+_r)}\frac{\partial G}{\partial\nu}\,d\s.
\]
From Lemma \ref{l:G} we know that $G$ is of the form
\[
G(x)=|x|^{-1}+{\cal R}(x)
\]
where ${\cal R}$ is regular. Since
\[
\int_{\G_2(B^+_r)}\frac{\partial}{\partial\nu}|x|^{-1}\,d\s=-\frac 12 \,|{\mathbb S}^{n-1}|
\]
and 
\[
\int_{\G_2(B_r^+)}\frac{\partial{\cal R}}{\partial\nu}\,d\s\ \mathop{\longrightarrow}\limits_{r\to+0^+}\ 0,
\]
we conclude that 
\[
\lim_{r\to0^+}\int_{\G_2(B^+_r)}\frac{\partial
  u}{\partial\nu}\,d\s=-\frac 12\,b\,|{\mathbb S}^{n-1}|
\]
thus getting the conclusion.
\end{pf}

\begin{center}

{\bf Acknowledgements}

\end{center}

The authors would like to thank the referee for pointing out a mistake in a previous version of this work.
They would like to thank Prof. A. Ambrosetti and Prof. Alice Chang for their interest in this work and for their constant support. V. F. is supported by M.U.R.S.T. under the national project ``Variational Methods and Nonlinear Differential Equations'' .

\vskip3truecm

\noindent {\sc Veronica Felli}\\
Scuola Internazionale Superiore di Studi Avanzati (S.I.S.S.A.)\\
Via Beirut 2-4, 34014 Trieste, Italy.\\
{\it E-mail: felli@sissa.it}.\\ \\
{\sc Mohameden Ould Ahmedou}\\
Rheinische Friedrich-Wilhelms-Universit\"{a}t Bonn\\
Mathematisches Institut\\
Beringstrasse 4 , D-53115 Bonn, Germany.\\
{\it E-mail: ahmedou@math.uni-bonn.de}.

\end{document}